\documentclass[12pt]{article}
\usepackage{verbatim}
\usepackage{latexsym,amsfonts}
\usepackage{amssymb,amsmath,amsthm}

\newtheorem{lemma}{Lemma}[section]
\newtheorem{thm}[lemma]{Theorem}
\newtheorem{cor}[lemma]{Corollary}
\newcommand{\F}{{\mathbb F}}
\newcommand{\Z}{{\mathbb Z}}

\title{Nonexistence of permutation binomials \\ of certain shapes}

\author{Ariane M. Masuda\thanks{This work proves the conjectures stated
in the first author's talk at the November 2006 BIRS workshop on
Polynomials over Finite Fields and Applications. The authors thank
BIRS for providing wonderful facilities.  The first author
was at Carleton University when this research was performed.} \\
\small{Department of Mathematics and Statistics}\\[-0.8ex]
\small University of Ottawa, Ottawa, ON K1N 6N5, Canada \\[-0.8ex]
\small \texttt{amasuda@uottawa.ca} \\ \\
Michael E. Zieve\footnotemark[1] \\
\small Center for Communications Research\\[-0.8ex]
\small 805 Bunn Drive, Princeton, NJ 08540 \\[-0.8ex]
\small \texttt{zieve@math.rutgers.edu}}

\date{}

\begin{document}
\maketitle

\begin{abstract}
Suppose $x^m+ax^n$ is a permutation polynomial over $\F_p$,
where $p>5$ is prime and $m>n>0$ and $a\in\F_p^*$.  We prove that
$\gcd(m-n,p-1)\notin\{2,4\}$.  In the special case that either
$(p-1)/2$ or $(p-1)/4$ is prime, this was conjectured in a recent
paper by Masuda, Panario and Wang.
\end{abstract}

\maketitle


\section{Introduction}

A polynomial over a finite field is called a \emph{permutation
polynomial} if it permutes the elements of the field.  These
polynomials have been studied intensively in the past two centuries.
Permutation monomials are completely understood: for $m>0$, $x^m$
permutes $\F_q$ if and only if $\gcd(m,q-1)=1$.  However, even
though dozens of papers have been written about them, permutation
binomials remain mysterious.  In this note we prove the following
result:
\begin{thm}
\label{main} If $p>5$ is prime and $f:=x^m+ax^n$ permutes $\F_p$,
where $m>n>0$ and $a\in\F_p^*$, then $\gcd(m-n,p-1)\notin\{2,4\}$.
\end{thm}
In case $(p-1)/2$ or $(p-1)/4$ is prime, this was conjectured in the
recent paper \cite{MPW} by Panario, Wang and the first author. It is
well-known that the gcd is not 1: for in that case, $f$ has more
than one root in $\F_p$, since $x^{m-n}$ is a permutation
polynomial. It is much more difficult to show that the gcd is not 2
or 4.

In Section~\ref{sec clique} we prove some general results about
permutation binomials, and in particular we show that it suffices to
prove Theorem~\ref{main} when $m-n$ divides $p-1$.  Then we prove
Theorem~\ref{main} in Section~\ref{sec main}.

Throughout this paper, we want to ignore permutation binomials that
are really monomials in disguise.  Here one can disguise a permutation
monomial (over $\F_q$) by adding a constant plus a multiple of $x^q-x$;
such addition does not affect the permutation property.
Thus, we say a permutation binomial
of $\F_q$ is \emph{trivial} if it is congruent modulo $x^q-x$ to the
sum of a constant and a monomial.  In other words, the nontrivial
permutation binomials are those whose terms have degrees being positive
and incongruent modulo $q-1$.


\section{Permutation binomials in general}
\label{sec clique}

\begin{lemma}\label{red1}
If $f$ is a permutation polynomial over\/ $\F_q$, then the greatest
common divisor of the degrees of the terms of $f$ is coprime to $q-1$.
\end{lemma}
\begin{proof}
Otherwise $f$ is a polynomial in $x^d$, where $d>1$ divides $q-1$,
but $x^d$ is not a permutation polynomial so $f$ is not one either.
\end{proof}

\begin{lemma}\label{clique}
Let $d\mid (q-1)$, and suppose there are no nontrivial permutation binomials
over\/ $\F_q$ of the form $x^e(x^d+a)$.
Then there are no nontrivial permutation binomials over\/ $\F_q$ of the form
$x^n(x^k+a)$ with $\gcd(k,q-1)=d$.
\end{lemma}
\begin{proof}
Suppose $f(x):=x^n(x^k+a)$ permutes $\F_q$, where $n,k,a\ne 0$. Let
$d=\gcd(k,q-1)$.  Pick $r>0$ such that $kr\equiv d\pmod{q-1}$ and
$\gcd(r,q-1)=1$. Then $f(x^r)$ permutes $\F_q$ and $f(x^r)\equiv
x^{nr}(x^d+a)\pmod{x^q-x}$.
\end{proof}

Lemma~\ref{clique} immediately implies the following
result from~\cite{MPW}:

\begin{cor}
If $q-1$ is a Mersenne prime, then there are no nontrivial permutation
binomials over\/ $\F_q$.
\end{cor}

We give one further reduction along the lines of Lemma~\ref{clique}:

\begin{lemma}\label{red3}
Let $d,n,e>0$ satisfy $d|(q-1)$, $\gcd(ne,d)=1$ and $n\equiv
e\pmod{(q-1)/d}$. Then $x^n(x^d +a)$ permutes\/ $\F_q$ if and only if
$x^e(x^d+a)$ does.
\end{lemma}
\begin{proof}
Write $f:=x^n(x^d+a)$ and $g:=x^e(x^d+a)$.  For any $z\in\F_q$ with
$z^d=1$, we have $f(zx)=z^n f(x)$; since $\gcd(n,d)=1$, this implies
that the values of $f$ on $\F_q$ comprise all the $d^{\operatorname{th}}$
roots of the values of $f(x)^d$.  Since $f(x)^d\equiv
g(x)^d\pmod{x^q-x}$, the result follows.
\end{proof}

Finally, since we constantly use it, we give here a version of
Hermite's criterion~\cite{H}:

\begin{lemma}
A polynomial $f\in\F_q[x]$ is a permutation polynomial if and only if
\begin{enumerate}
\item for each $i$ with $0<i<q-1$, the reduction of
$f^i$ modulo $x^q-x$ has degree less than $q-1$; and
\item $f$ has precisely one root in $\F_q$.
\end{enumerate}
\end{lemma}


\section{Proof of Theorem~\ref{main}}
\label{sec main}

In this section we prove Theorem~\ref{main}.  We treat the cases of gcd 2
and 4 separately.

\begin{thm}
\label{2p+1} If $p$ is prime and $x^n(x^k+a)$ is a nontrivial permutation
binomial over\/ $\F_p$, then $\gcd(k,p-1)>2$.
\end{thm}
\begin{proof}
There are no nontrivial permutation binomials over $\F_2$ or $\F_3$,
so we may assume $p=2\ell+1$ with $\ell>1$.  By Lemma~\ref{clique},
it suffices to show there are no nontrivial permutation binomials of the form
$f:=x^n(x^d+a)$ with $d\in\{1,2\}$. This is clear for $d=1$ (since then
$f(0)=f(-a)$), so we need only consider $d=2$. Assume
$f:=x^n(x^2+a)$ is a permutation binomial. Lemma~\ref{red1} implies
$n$ is odd.

Suppose $\ell$ is odd.  We will use
Hermite's criterion with exponent $\ell-1$; to this end, we compute
\[
f^{\ell-1} = x^{n\ell-n}  (x^2 + a)^{\ell-1}= x^{n\ell-n} \sum_{i=0}^{\ell-1}
\binom{\ell-1}{i} a^{\ell-1-i} x^{2i}.
\]
\noindent Write $f^{\ell-1} =
\sum_{i=0}^{\ell-1} b_ix^{n\ell-n+2i}$, where $b_i = \binom{\ell-1}{i}
a^{\ell-1-i}$.  Since $\ell-1<p$ and $p$ is prime, each $b_i$ is nonzero.
Thus, the degrees of the terms of $f^{\ell-1}$ are precisely the elements
of
\[
S=\{n\ell-n,n\ell-n+2,n\ell-n+4,\ldots,n\ell-n+2\ell-2\}.
\]
\noindent Since  $\ell$ is odd,  $S$ consists of $\ell$ consecutive
even numbers, so it contains a unique multiple of $p-1=2\ell$. Thus
the reduction of $f^{\ell-1}$ modulo $x^p-x$ has degree $p-1$, which
contradicts Hermite's criterion.

If $\ell$ is even then $f^{\ell}=\sum_{i=0}^{\ell}c_ix^{n\ell+2i}$, where each
$c_i=\binom{\ell}{i} a^{\ell-i}$ is nonzero.  The degrees of the terms of
$f^{\ell}$ consist of the $\ell+1$ consecutive even numbers $n\ell, n\ell+2,
\ldots, n\ell+2\ell$. Since $n$ is odd, $n\ell$ is not a multiple of
$p-1=2\ell$.  Thus $f^{\ell}$ has a unique term of degree divisible by $p-1$,
which again contradicts Hermite's criterion.
\end{proof}

\begin{thm}
\label{4p+1} If $p$ is prime and $x^n(x^k+a)$ is a nontrivial
permutation binomial over\/ $\F_p$, then $\gcd(k,p-1)\neq  4$.
\end{thm}

\begin{proof}
Plainly we need only consider primes $p$ with $p\equiv 1\pmod{4}$.
By Lemma~\ref{clique}, it suffices to show there are no nontrivial
permutation binomials of the form $x^n(x^4+a)$.  By Lemma~\ref{red1},
we may assume $n$ is odd.  By Lemma~\ref{red3}, it suffices to show
nonexistence with $0<n< (p-1)/4$ if $p\equiv 1\pmod{8}$, and with
$0<n< (p-1)/2$ if $p\equiv 5\pmod{8}$.  Assume $f:=x^n(x^4+a)$ is a
nontrivial permutation binomial with $n$ satisfying these
constraints.

First suppose $p\equiv 1\pmod8$, say $p=8\ell+1$; here our assumption is
$0<n<2\ell$.  The set of degrees of terms of $f^{2\ell}$ is
\[
S=\{2\ell n,2\ell n+4,2\ell n+8,\ldots,2\ell n+8\ell\}.
\]
\noindent When $\ell$ is even, $S$ consists of $2\ell+1$ consecutive multiples
of $4$.  Since $n$ is odd, $2\ell n$ is not a multiple of $8\ell$, so $S$
contains precisely one multiple of $p-1=8\ell$, contradicting Hermite's
criterion.  So assume $\ell$ is odd; since $8\ell+1$ is prime, we have
$\ell\ge 5$.  Now the set of degrees of terms of $f^{2\ell+2}$ is
\[
S=\{2\ell n+2n,2\ell n+2n+4,2\ell n+2n+8,\ldots,2\ell n+2n+4(2\ell+2)\}.
\]
Here $S$ consists of $2\ell+3$ consecutive multiples of $4$, so it contains
a multiple of $p-1=8\ell$. By Hermite's criterion, $S$ must have at least
two such multiples.  Thus, $8\ell$ divides either $2\ell n+2n$, $2\ell n
+2n+4$ or $2\ell n+2n+8$, so $\ell$ divides either $n$, $n+2$ or $n+4$.
Since $\ell\ge 5$ and $0<n<2\ell$, we have $n+4<3\ell$; since $n$ is odd,
it follows that $\ell$ equals either $n$, $n+2$ or $n+4$.  But then
$f^8$ has a unique term of degree divisible by $p-1=8\ell$, contradicting
Hermite's criterion.

Thus we have $p\equiv 5\pmod{8}$; write $p=4\ell+1$ with $\ell$ odd, where
again $0<n<2\ell$.
Suppose $\ell\equiv 1\pmod{4}$. If $\ell=1$ then $f$ is trivial, so assume
$\ell>1$.  The set of degrees of terms of $f^{\ell-1}$ is
\[
S=\{n\ell-n,n\ell-n+4,n\ell-n+8,\ldots,n\ell-n+4\ell-4\}.
\]
\noindent Since  $\ell\equiv 1 \pmod 4$, the set $S$ consists of $\ell$
consecutive multiples of $4$, so $S$ contains precisely one multiple
of $p-1=4\ell$, contradicting Hermite's criterion.

Thus $\ell\equiv 3\pmod{4}$. The set of degrees of terms of
$f^{\ell+1}$ is
\[
S=\{n\ell+n,n\ell+n+4,n\ell+n+8,\ldots,n\ell+n+4\ell+4\}.
\]
\noindent Since $S$ consists of $\ell+2$ consecutive multiples of
$4$, it certainly contains a multiple of $4\ell$, so (by Hermite's
criterion) it must contain two such multiples. Thus either
$n(\ell+1)$ or $n(\ell+1)+4$ is a multiple of $4\ell$, so $\ell$
divides either $n$ or $n+4$. Since $n$ is odd and $0<n< 2\ell$, the
only possibilities are $n=\ell$ or $n=\ell-4$ or $(n,\ell)=(5,3)$.
If $n=\ell-4$ then $f^4$ has degree $4\ell=p-1$, contradicting
Hermite's criterion.  If $(n,\ell)=(5,3)$, then $p=13$ and
$a^{-1}f(x^{11})$ permutes $\F_p$; since $a^{-1}f(x^{11})\equiv
x^3(x^4+a^{-1})\pmod{x^{13}-x}$, it suffices to treat the case
$n=\ell$. Finally, suppose $n=\ell$, so $f=x^{\ell}(x^4+a)$ permutes
$\F_p$. The degrees of the terms of $f^4$ are
\[
4\ell, 4\ell+4, 4\ell+8, 4\ell+12, 4\ell+16.
\]
\noindent We have our usual contradiction if the degree $4\ell$ term
is the unique term of $f^4$ with degree divisible by $4\ell$, so the
only remaining possibility is that $4\ell$ divides either $4$, $8$,
$12$ or $16$.  Since $\ell\equiv 3\pmod{4}$, the only possibility is
$\ell=3$.  Finally, when $\ell=3$, the coefficient of $x^{12}$ in
the reduction of $f^4$ modulo $x^{13}-x$ is $a^4+4a$, which must be
zero (by Hermite), so $a^3=-4$; but the cubes in $\F_{13}^*$ are
$\pm 1$ and $\pm 8$, contradiction.
\end{proof}



\end{document}